\documentclass{amsart}
\newtheorem{theorem}{Theorem}
\newtheorem{lemma}{Lemma}
\newtheorem{proposition}{Proposition}
\theoremstyle{definition}
\newtheorem{definition}{Definition}

\newcommand{\wg}{\widehat{G}}
\newcommand{\su}[1]{\mathrm{Supp}\,#1}

\newcommand{\bz}{\bar{0}}
\newcommand{\bo}{\bar{1}}
\theoremstyle{remark}
\newtheorem{remark}{Remark}
\newcommand{\z}{\mathbb{Z}_2}
\newcommand{\vp}{\varphi}
\newcommand{\aut}[1]{\mathrm{Aut}\,(#1)}
\newcommand{\baut}[1]{\overline{\mathrm{Aut}}\,(#1)}
\begin{document}
\title{Group gradings on superinvolution simple superalgebras}
\author{Yu. Bahturin, M.Tvalavadze, T. Tvalavadze}
\address{Department of Mathematics and Statistics,
Memorial University of Newfoundland, St.John's, NL, A1C5S7,Canada}
\email{yuri@math.mun.ca}
\thanks{This work was supported in part by NSERC grant
227060-04}
\subjclass[2000]{Primary 16W10, 16W50; Secondary 16W55}
\date{}

\keywords{Associative superalgebras, superinvolution, gradings}

\begin{abstract}
In this paper we describe all group gradings by an arbitrary
finite group $G$ on non-simple finite-dimensional superinvolution
simple associative superalgebras over an  algebraically closed
field $F$ of characteristic 0 or coprime to the order of $G$.
\end{abstract}

\maketitle

\section*{Introduction}

In the paper \cite{bg}, Yu.Bahturin and A. Giambruno described the
group gradings by finite abelian groups $G$  on the matrix
algebra $M_n(F)$ over an algebraically closed field $F$ of
characteristic different from 2, which are  respected by an
involution. Besides, under some restrictions on the base field,
they classified all $G$-gradings on all finite-dimensional
involution simple algebras.

In this paper we deal with finite-dimensional  associative
superalgebras that are simple with respect to some superinvolution
$*$ over an algebraically closed field of characteristic zero or
coprime to the order of $G$. First, we give a description of such
associative superalgebras. Second, we classify all group gradings
on $*$-simple associative superalgebras that are not simple
associative algebras.

In the same way as the description of involution gradings on
involution simple associative algebras is important for the
determination of group gradings on classical simple Jordan and Lie
algebras  \cite{BShZ, ant} the description of superinvolution
gradings on superinvolution simple associative algebras is
important for the determination of groups gradings on simple
Jordan and Lie superalgebras. The case of superalgebras that are simple
algebras is due to the second and the third authors, and is to be submitted for publication shortly.

\section{Definitions and introductory remarks}

Let $R$ be an associative superalgebra, or, in other words, an associative algebra
with a fixed  $\mathbb{Z}_2$-grading $R=R_{\bar 0}\oplus R_{\bar 1}$. Since all
algebras and superalgebras considered in this paper are associative we will normally
drop the word associative in what follows. Also, if not stated otherwise, all subalgebras,
ideals and homomorphisms are $\z$-graded. We say that $R$ is simple if it has no nontrivial proper ($\z$-graded) ideals.
It is well-known \cite{W} that any finite-dimensional simple (associative) superalgebra over
an algebraically closed field of characteristic different from 2 is isomorphic to either
$M_{k,l}(F)$, the full matrix algebra $M_n(F)$ with a $\z$-grading completely
determined by two nonnegative integers $k,l$, $k+l=n$, or a subalgebra $R=Q(n)$ of $M_{2n}(F)$ consisting of all matrices of the form
$\left(\begin{array}{cc}
             X& Y\\
             Y& X
             \end{array}\right)$ with $R_{\bar 0}=\left(\begin{array}{cc}
             X& 0\\
             0& X
             \end{array}\right)$ and $R_{\bar 1}=\left(\begin{array}{cc}
             0& Y\\
             Y& 0
             \end{array}\right)$. A convenient notation for $R=Q(n)$
            is $R=A\oplus tA$ where $A\cong M_n(F)$,
$t^2=1$. Then $R_{\bar 0}=A$ and $R_{\bar 1}=tA$.

\begin{definition} Let $R$ be a superalgebra. A
\emph{superinvolution} on $R$ is a $\mathbb{Z}_2$-graded linear map
$*:R\to R$ such that $(x^*)^*=x$ for all $x\in R$ and
$(xy)^*=(-1)^{|x||y|}y^*x^*$ for all homogeneous $x,y\in R$, of degrees $|x|$ and $|y|$, respectively. A more general notion is that of \emph{superantiautomorphism}, that is, a linear map $\vp:R\rightarrow R$ such that $\vp(xy)=(-1)^{|x||y|}\vp(y)\vp(x)$ for all homogeneous $x,y\in R$, as above.
\end{definition}

In this paper we will be interested in superinvolution simple superalgebras.

\begin{definition}
Let $(R,*)$ be an superalgebra endowed with a superinvolution
$*$. We say that $R$ is \emph{superinvolution simple} if $R^2\neq \{ 0\}$ and $R$ has no $\z$-graded ideals stable under $*$.
\end{definition}

If $*$ is a superinvolution on $R$,
then obviously the restriction of $*$ to $R_{\bz}$ is an ordinary involution. The same is true for superantiautomorphisms.
Thus, there is no confusion to abbreviate the terms superinvolution and superantiautomorphism to involution and antiautomorphism, respectively. In particular, superinvolution simple superalgebras will be simply called involution simple.

The following are examples of involution simple superalgebras.

\par\medskip
{\large Examples} {\bf 1.} The {\it orthosymplectic}
involution on $R=M_{r,2s}(F)$ is given by
$$ \left(\begin{array}{cc}
             X& Y\\
             Z& T
             \end{array}\right)^{osp}=
   \left(\begin{array}{cc}
             I_r& 0\\
             0&Q
             \end{array}\right)^{-1}
    \left(\begin{array}{cc}
             X& -Y\\
             Z&  T
             \end{array}\right)^t
     \left(\begin{array}{cc}
             I_r& 0\\
             0& Q
             \end{array}\right)
$$
where $t$ denotes the usual matrix transpose,
 $Q=\left(\begin{array}{cc}
           0& I_s\\
        -I_s& 0
     \end{array}\right)$, and $I_r, I_s$ are
the identity matrices of orders $r$,$s$, respectively.
\par\medskip
{\bf 2.} Let us consider $R=M_{r,r}(F)$. We will call the following
involution defined on $M_{r,r}(F)$ the {\it transpose}
involution:
$$  \left(\begin{array}{cc}
            X& Y\\
            Z& T
            \end{array}\right)^{trp}=
            \left(\begin{array}{cc}
            T^t& -Y^t\\
            Z^t& X^t
            \end{array}\right).
$$
\par\medskip
{\bf 3.} Let $A$ be a superalgebra. Consider a new
 superalgebra $A^{sop}$ which has the same
$\mathbb
{Z}_2$-graded vector space structure  as $A$  but the product
of $A^{sop}$ is given on homogeneous $a,b$ of degrees $|a|$, $|b|$ by
$$ a\circ b =(-1)^{|a||b|}ba.$$
%where $|a|=\left\{\begin{array}{cc}
%                 0,&\, if\, a\in A_0\\
%                 1,&\, if\, a\in A_1
%            \right.$

Let $S=A\oplus A^{sop}$ be the direct sum of two ideals $A$ and
$A^{sop}$. This is a $\mathbb{Z}_2$-graded algebra with  $R_{\bar
0}=A_{\bz}\oplus A_{\bz}^{sop}$,  $R_{\bar 1}=A_{\bo}\oplus A_{\bo}^{sop}$. We
denote an arbitrary element $x$ from $R$ as a pair of elements
from $A$, i.e. $x=(a,b)$ where $a,b\in A$. The product in $R$ is
given by
$$ (a_0+a_1,b_0+b_1)\cdot
(a'_0+a'_1,b'_0+b'_1)=(a_0a'_0+a_1a'_1+a_0a'_1+a_1a'_0,b'_0b_0-b'_1b_1+b'_1b_0+b'_0b_1),$$
where $a_0,b_0,a'_0,b'_0\in A_0$, $a_1,b_1,a'_1,b'_1\in A_{\bo}$.

A linear mapping defined by $(a,b)^{ex}=(b,a)$ is an involution called {\it
exchange involution}. If $A$ is simple then $(S,ex)$ is an involution simple superalgebra.

\begin{definition}
Let $R$ and $S$ be two superalgebras endowed with involutions $\ast$ and $\dagger$. We say that $(R,\ast)$ and $(S,\dagger)$ are \emph{isomorphic} if  there exists an isomorphism of superalgebras $\vp:R\rightarrow S$ such that $\vp(x^\ast)=\vp(x)^\dagger$ for all $x\in R$. If $R=S$ then $\vp$ is an automorphism of $R$ and $\ast$, $\dagger$ are called \emph{conjugate} by $\vp$. In this case we have
$\dagger=\varphi\circ\ast\circ\varphi^{-1}$.
\end{definition}

If $(R,\ast)$ is an involution simple superalgebra then a standard argument shows that either $R$ is a simple superalgebra or else there is a ($\z$-graded) ideal $A$ in $R$ such that $R=A\oplus A^\ast$. In the latter case the mapping $\vp: R\rightarrow S$ defined by $\vp(a+b^\ast)=(a,b)$ where $a,b\in A$ defines an isomorphism of involution simple superalgebras between $(R,\ast)$ and a standard superalgebra $(S,ex)$ of Example 3 above, where $A$ is simple.

 In \cite{r} M.Racine described all types of involutions on
 $A=M_{n,m}(F)=A_{\bz}+A_{\bo}$. It appears that if $\varphi$ is an
 involution on $A$ such that $A_{\bz}$ is an involution simple
 algebra under $\varphi$ restricted to $A_{\bz}$, then $n=m$
 and $\varphi$ is conjugate to the transpose involution.
 Otherwise, $\varphi$ is conjugate to the orthosymplectic
 involution. Also, it was shown in \cite{gs}, a superalgebra of the type
$Q(n)$ has no involutions.
We can summarize all the remarks above as the following.

\begin{proposition}\label{p1}
Any finite-dimensional involution simple superalgebra over an algebraically closed field of characteristic different from 2 is isomorphic
to one of the following:
\begin{enumerate}
 \item $R=M_{n,m}(F)$ with the orthosymplectic or transpose
involution.

    \item  $R=M_{n,m}(F)\oplus M_{n,m}(F)^{sop}$ with the ordinary
exchange involution.

    \item  $R= Q(n)\oplus Q(n)^{sop} $ with the ordinary exchange
involution.
\end{enumerate}
\end{proposition}

\section{Group Gradings}

One can define gradings of superagebras by the elements of very general sets with operations but as it turns out if the superalgebra is involution simple we can restrict ourselves to the case of abelian groups. A phenomenon of this kind was, probably, first mentioned in  \cite{PZ}. In the case of involutions see \cite{BShZ} and \cite{bg}.

\begin{definition} Given a semigroup $G$ and a superalgebra $R$ we say that $R$ is \emph{graded} by $G$ if $R=\bigoplus_{g\in G}R_g$ where each $R_g$ is a $\z$-graded vector subspace and $R_gR_h\subset R_{gh}$, for any $g,h\in G$. The subset $\su{R}=\{g\in G\,|\, R_g\neq\{0\}\}$ is called the \emph{support} of the grading.
\end{definition}

A semigroup with 1 is called \emph{cancellative} if each of $xg=xh$, $gx=hx$ implies $g=h$, for any $x,g,h\in G$.

\begin{proposition}
Let $R$ be a $G$-graded superalgebra, $G$ a cancellative semigroup. Suppose $R$ has an involution $*$
compatible with  this grading, that is, $R^*_g=R_g$, for any $g\in
G$, and also that $R$ is $*$-simple. Then, given any $g,h\in
\mbox{Supp}\,R$ we have that $gh=hg$. If, additionally, $1\in \su{R}$ then any $g\in \su{R}$ is invertible.
\end{proposition}

\begin{proof}
Let $g,h\in\mbox{Supp}\ R$. Suppose  $R_g R_h\ne 0$. We have to
show that $(R_g R_h)^*\subseteq R_{hg}$. Since  we deal with a superalgebra $G$-grading, $R_g=R_g^{\bz}+R_g^{\bo}$ and
$R_h=R_h^{\bz}+R_h^{\bo}$ where $R_g^{\bz}$, $R_h^{\bz}$ are even components,
$R_g^{\bo}$, $R_h^{\bo}$ are odd components. It follows from $
R_gR_h=(R_g^{\bz}+R_g^{\bo})(R_h^{\bz}+R_h^{\bo})\subseteq
R_g^{\bz}R_h^{\bz}+R_g^{\bo}R_h^{\bo}+R_g^{\bz}R_h^{\bo}+R_g^{\bo}R_h^{\bz}$ that
$(R_gR_h)^*\subseteq
(R_h^{\bz})^*(R_g^{\bz})^*+(R_h^{\bo})^*(R_g^{\bo})^*+(R_h^{\bo})^*(R_g^{\bz})^*+(R_h^{\bz})^*(R_g^{\bo})^*=
R_h^{\bz}R_g^{\bz}+R_h^{\bo}R_g^{\bo}+R_h^{\bo}R_g^{\bz}+R_h^{\bz}R_g^{\bo}\subseteq R_{hg}$. On
the other hand, $R_gR_h\subseteq R_{gh}$, $(R_gR_h)^*\subseteq
R_{gh}^*=R_{gh}$. Hence, $R_{gh}=R_{hg}$, $gh=hg$.

Now, pick $g,h\in \mbox{Supp}\,R$, and consider
$I=R_g+RR_g+R_gR+RR_gR$. It is easily seen that $I$ is a graded
ideal. Next we want to show that $I^*=I$. Since $RR_g=\sum_l
R_lR_g$, $(\sum_l R_gR_l)^*\subseteq \sum_l
R_l^{\bz}R_g^{\bz}+R_l^{\bo}R_g^{\bo}+R_l^{\bo}R_g^{\bz}+R_l^{\bz}R_g^{\bo}=\sum_l
(R_l^{\bz}+R_l^{\bo})(R_g^{\bz}+R_g^{\bo})=\sum_l R_lR_g=RR_g$. In a similar
manner we can show that
$(RR_gR)^*=(\sum_{l,k}R_lR_gR_k)^*=\sum_{l,k}
(R_lR_gR_k)^*=\sum_{k,l} R_kR_gR_l$. Therefore, $I$ is a graded
$*$-invariant non-zero ideal, hence, $I=R$. In particular,
$R_h\subset R_g+RR_g+R_gR+RR_gR$. The homogeneous components on
the right-hand side are of one of the forms: $g$, $kg$, $gl$,
$pgq$, for some $k,l,p,q\in G$. So, $h$ is one of these forms. It
follows that one of the spaces $R_g$ (if $g=h$), or $R_kR_g$, or
$R_gR_l$, or $R_pR_gR_g$ is different from zero, with either
$h=g$, or $h=kg$, or $h=gl$, or $h=pgq$. The case $h=g$ being
trivial, if $R_kR_g\ne 0$ with $h=kg$ then $kg=gk$ by what  was
proven before and then $hg=(kg)g=g(kg)=gh$, as needed. Similarly,
if $R_gR_l\ne 0$ with $gl\ne 0$. Now if $R_pR_gR_q\ne 0$ with
$h=pgq$, then $R_pR_g\ne 0$ and $R_gR_q\ne 0$ so that $pg=gp$ and
$gq=qg$. Again, $hg=(pgq)g=(pg)(qg)=gpgq=gh$, as required.

\medskip

The invertibility claim follows in exactly the same way as in \cite[Proposition 1]{semigr}.

\end{proof}

As a result, using Proposition \ref{p1}, we will assume in what follows, that we deal with abelian group gradings of finite-dimensional involution simple superalgebras. Actually, we restrict ourselves to the case where $G$ is finite and $R$ is not simple as a superalgebras (Cases (2) and (3) of Proposition \ref{p1}). As mentioned earlier, Case (1) is to be published in the joint paper of the second and the third authors.

\begin{remark}\label{r1} If $A$ is a superalgebra graded by an \emph{abelian} group $G$ then
the same homogeneous subspaces $A_g$, $g\in G$, define in
$A^{sop}$ a $G$-grading. We will denote these subspaces by
$A^{sop}_g$.
\end{remark}

The techniques we are going to use impose a further restriction on
the ground field $F$. Namely, we are going to use the
correspondence between the gradings on a (super) algebra $R$ by a
finite abelian group $G$ and the actions on $R$ of the dual group
$\wg$ by automorphisms (see, for example, \cite{bsz}). For this to
work, we need to make sure that if the order of $G$ is $d$ then
$F$ contains $d$ different roots of $1$ of degree $d$. If this
condition holds then each grading $R=\bigoplus_{g\in G}R_g$
defines a homomorphism $\alpha:\wg\rightarrow\mathrm{Aut}\,R$
given by $\alpha(\chi)(r)=\chi(g)r$ provided that $r\in R_g$,
$g\in G$. Also the grading can be recovered if we have a
homomorphism $\alpha$, as above.

We start with a general result (Exchange Theorem below)
obtained by the first author. An important particular case can be found in \cite{bs}. Let $G$ be a finite abelian group and $V$
a vector space. Suppose we have two $G$-gradings on $V$:
$$ V=\oplus_{g\in G}\overline{V}_g,\quad \alpha:\wg\to\mbox{Aut}\,V,\eqno (1')$$
$$ V=\oplus_{g\in G}\widetilde{V}_g,\quad \beta:\wg\to\mbox{Aut}\,V,\eqno (2')$$

\noindent where $\alpha,\beta:\wg\to\mbox{Aut}\,V$ are
homomorphisms of the dual group $\wg$ corresponding to the
above gradings in the following way. Given $\chi\in\wg$ we
define $\alpha(\chi)$ on an element $v$ of $\overline{V}_g$, for each
$g$, by $\alpha(\chi)(v)=\chi(g)v$. Similarly for $(2')$. Suppose
$\Lambda\subset \wg$ is a subgroup such that
$\alpha(\lambda)=\beta(\lambda)$, for each $\lambda\in \Lambda$.
Let us denote by $H$ the orthogonal complement
$\Lambda^{\bot}=\{g\in G|\,\lambda(g)=1,\,\lambda\in\Lambda\}$.
Assume further that the subgroups $\alpha(\wg)$ and $\beta(\hat
G)$ commute elementwise.

Let us consider a homomorphism $\gamma:\wg\to \mbox{Aut}\,V$
given by $\gamma(\chi)=\alpha^{-1}(\chi)\beta(\chi)$. In this case
we can define $H$-grading of $V$ as follows:
$V^{(h)}=\{v|\,\gamma(\chi)(v)=\chi(h)v, \chi\in\wg\}$.
\par\medskip
\noindent{\bf Theorem} (Exchange Theorem). {\it The three
gradings defined above are connected by the following equations
$$ \overline{V}_g=\oplus_{h\in H}(\widetilde{V}_{gh}\cap V^{(h)}),\quad
\widetilde{V}_g=\oplus_{h\in H}(\overline{V}_{gh}\cap V^{(h^{-1})})\eqno
(3')$$ If $V$ is an algebra and $(1')$, $(2')$ are algebra
gradings, then $(3')$ are relations for the algebra gradings.}
\begin{proof}
Let us prove the first equality. Since all gradings are
compatible, we have $\overline{V}_g=\oplus_{h\in H}(\overline{V}_{g}\cap
V^{(h)})$. Thus it is enough to prove, for any $g\in G, h\in H$,
that $\widetilde{V}_{gh}\cap V^{(h)} =\overline{V}_{g}\cap V^{(h)}.$ If $v\in
\widetilde{V}_{gh}\cap V^{(h)}$ then $\beta(\chi)(v)=\chi(gh)v$ and
$\gamma(\chi)(v)=\chi(h)v$. Hence also
$\gamma(\chi)^{-1}(v)=\chi(h)^{-1}v$. Now
$$
\alpha(\chi)(v)=\alpha(\chi)\beta(\chi)^{-1}\beta(\chi)(v)=
\gamma(\chi)^{-1}\beta(\chi)(v)=\chi(h)^{-1}\chi(gh)v
$$
 proving
$\widetilde{V}_{gh}\cap V^{(h)}\subset \overline{V}_g\cap V^{(h)}$.

If $b\in \overline{V}_g \cap V^{(h)}$ then $\alpha(\chi)(b)=\chi(g)b, \
\gamma(\chi)(b)=\chi(h)b$. Therefore
$$
\beta(\chi)(b)=\alpha(\chi)\alpha(\chi)^{-1}\beta(\chi)(a)=
\alpha(\chi)\gamma(\chi)(a) =\chi(g)\chi(h)a=\chi(gh)a.
$$
It follows that $\overline{V}_g \cap V^{(h)}\subset \widetilde{V}_{gh} \cap
V^{(h)}$. Finally, $\overline{V}_g \cap V^{(h)}= \widetilde{V}_{gh} \cap
V^{(h)}$ for any $g\in G$ and thus we have the first equality in
$(3')$. The second is similar. It is easy to check that if $V$ is
an algebra and ($1'$), ($2'$) are algebra gradings, then ($3'$)
provides us with the relations between algebra gradings as well.
The proof is complete.
\end{proof}

One of important tools in the proof of the main results of our work is a recent result from \cite{ant}, as follows.

\begin{theorem}\label{tgrant} Let $M_n(F)=A=\bigoplus_{g\in G}A_g$ be a $G$-grading on $M_n(F)$ over a field $F$
of characteristic not 2, which contains $d$ different roots of 1, $d=|G|$. Suppose there is a graded antiautomorphism
$\vp$ whose restriction to $R_e$ is an involution. Then there is a $G$-graded automorphism $\psi$ of $R$ such that
$\vp\psi=\psi\vp$ and $\psi^2=\vp^2$.
\end{theorem}

A consequence of this result which interests us is as follows. Let us denote by $\baut{A}$ the group of automorphisms and antiautomorphisms of $A$. In the case $A=M_n(F)$, $[\baut{A}:\aut{A}]=2$.

\begin{theorem}\label{tautant} Let $P$ be a finite abelian subgroup in $\baut{A}$, $A=M_n(F)$ over a
field $F$ of characteristic not 2, which contains $d$ different roots of 1, $d=|G|$. S
uppose $\vp\in P\setminus \aut{A}$. Then there exists $\psi\in\aut{A}$
commuting with all elements in $P$ and $\psi^2=\vp^2$.
\end{theorem}

\begin{proof} Set $Q=P\cap\aut{A}$. Then $Q$ is a subgroup of index 2 in $P$. Let $G$ be a finite abelian group whose dual is $Q$. That is, the elements of $Q$ can be viewed as multiplicative characters on $G$. As noted earlier, in this case $A$ becomes $G$-graded if one sets $A_{g}=\{a\in A\,|\,\chi(a)=\chi(g)a\mbox{ for any }\chi\in Q\}$. Since $\vp$ commutes with the elements of $Q$, the antiautomorphism $\vp$ is a $G$-graded map. Also, because $\vp^2\in Q$, we have that the restriction of $\vp$ to $R_e$ is an involution. Applying Theorem \ref{tgrant}, we find a $G$-graded automorphism $\psi$ such that $\vp\psi=\psi\vp$ and $\psi^2=\vp^2$. Now if $\chi$ is an arbitrary element of $Q$ and $a$ a homogeneous element of degree $g\in G$ then $\psi(\chi(a))=\psi\chi(g)a=\chi(g)\psi(a)=\chi(\psi(a))$ because $\psi(a)\in A_g$. It follows that $\psi\chi=\chi\psi$ and $\psi$ commutes with all elements of $P$, as required.
\end{proof}

\section{Antiautomorphisms of graded superalgebras}\label{sags}

Theorem \ref{tgrant} is no longer true in the case of (super) antiautomorphisms of matrix superalgebras. The simplest example is the trivial grading and the (super) antiautomorphism defined on $A=M_{n,m}$, $n,m$ odd, by $$ \vp(X)=\left(\begin{array}{cc}
             A& -B\\
             C&  D
             \end{array}\right)^t.$$
             Luckily, the argument of \cite{ant} can be adapted to the case of superalgebras although we have to deal with higher powers of the antiautomorphisms in question. We start with a generalization of the results of \cite{BShZ} about fine involution gradings.

\begin{theorem}
Let $R=M_{n,m}(F)$, $n,m\ge 1$,  be a non-trivial matrix
superalgebra with an antiautomorphism  $\varphi$ over an
algebraically closed field $F$ of characteristic zero or coprime
to the order of $G$, where  $G$ is a finite Abelian group. Then
$R$ admits no fine $G$-gradings respected by $\varphi$.
\end{theorem}
\begin{proof}
Assume the contrary,  that is $R=\oplus_{g\in G} R_g$ is  a fine
$G$-grading respected by $\varphi$, $\varphi(R_g)=R_g$. Since $R$
is a superalgebra with a fine $G$-grading, according to
\cite{BSh},
 $\dim{R_{\bar{0}}}=
\dim{R_{\bar{1}}}$, that is, $n=m$. Let $R_{\bar 0}$ be denoted by
$A$. Then
$$ A=\oplus_{g\in G}\, A_g,$$
where $A_g=R_{\bar 0}\cap R_g$. This grading is also fine and
compatible with $\varphi$. Notice that $A=I_1\oplus I_2$, the sum
of two isomorphic simple ideals.

Next let $\widehat G$ be the dual group of $G$, and
$\alpha:{\widehat G}\to\mbox{Aut}\,A$ the homomorphism
accompanying our grading. If for each $\eta \in \widehat G$,
$\alpha(\eta)(I_i)=I_i$, then a fine $G$-grading of $A$ induce
$G$-grading on both ideals such that $A_g=(I_1)_g\oplus(I_2)_g$.
In particular, $A_e=(I_1)_e\oplus (I_2)_e$, $(I_i)_e\ne \{0\}$.
This contradicts the fact that our $G$-grading is fine. Therefore,
there exists $\xi\in \widehat G$ such that $\alpha(\xi)(I_1)=I_2$.
Hence, ${\widehat G}=\Lambda\cup \Lambda\xi$ where
$\Lambda=\{\eta\in {\widehat G}|\alpha(\eta)(I_i)=I_i\}$ and
$\xi^2\in \Lambda$. Then $H=\Lambda^{\bot}$ is a subgroup of $G$
of order 2 and $\widehat{G/H}\cong\Lambda$. Let $H=\{e,h\}$ where
$h^2=e$. Next we can consider the induced
$\overline{G}=G/H$-grading of $A$. Let $\bar g=gH$ for any $g\in
G$. Then $A_{\bar g}=A_g+A_{gh}$. Since $\widehat{G/H} * I_i =
\Lambda
* I_i =I_i$ where $i \in {1,2 }$, $I_i$ is a ${G/H}$-graded ideal.
It follows from $A_{\bar e}=(I_1)_{\bar e}\oplus (I_2)_{\bar e}$,
$(I_i)_{\bar e}\ne\{0\}$, and $\dim\,A_{\bar e}=2$ that
$\dim\,(I_i)_{\bar e}=1$. Therefore, both ${G/H}$-gradings on
$I_1$ and $I_2$ are fine.

Next we following two cases may occur.

{\it Case 1.} Let $\varphi(I_1)=I_2$. Notice that $A_{\bar
g}=(I_1)_{\bar g}\oplus (I_2)_{\bar g}$ for each $\bar g\in
\overline{G}$. Following arguments in \cite{bs}, we can recover
our original $G$-grading. In fact,
$$ A_g=\{X+\xi(g)^{-1}(\xi*X)|\,X\in A_{\bar g}\}.\eqno (1)$$
For example, let us take $ X=\left(\begin{array}{cc}
                                   X_{\bar g}& 0\\
                                        0    &  0
                                    \end{array}\right),$
$X_{\bar g}\in (I_1)_{\bar g}$. Then, by (1),
$$ 0\ne X+\xi(g)^{-1}(\xi*X)=\left(\begin{array}{cc}
                                   X_{\bar g}& 0\\
                                        0    &\xi(g)^{-1}(\xi*X_{\bar g})
                                    \end{array}\right)\in A_g.$$
Since $\dim\,A_g=1$, $A_g=\mbox{span}
\left\{\left(\begin{array}{cc}
                                   X_{\bar g}& 0\\
                                        0    &\xi(g)^{-1}(\xi*X_{\bar
                                        g})\end{array}\right)\right\}$.
Recall that $\varphi(A_g)=A_g$ where $\varphi$ can be represented
as follows:
$$ \varphi*\left(\begin{array}{cc}
                                   X& 0\\
                                   0 &  Y
                                    \end{array}\right)=
                                                 \left(\begin{array}{cc}
                                   \varphi_0(Y)& 0\\
                                   0 &  \varphi_1(X)
                                    \end{array}\right),$$
   where $\varphi_0$ and $\varphi_1$ are antiautomorphisms. Hence
   $$
                        \varphi* \left(\begin{array}{cc}
                                   X_{\bar g}& 0\\
                                        0    &\xi(g)^{-1}(\xi*X_{\bar
                                        g})\end{array}\right)=
                                  \left(\begin{array}{cc}
                     \xi(g)^{-1}(\varphi_0\xi)*(X_{\bar g})&
                     0\\
                     0& \varphi_1(X_{\bar g})
                     \end{array}\right)=
                     $$
                     $$\lambda_g \left(\begin{array}{cc}
                                   X_{\bar g}& 0\\
                                        0    &\xi(g)^{-1}(\xi*X_{\bar
                                        g})\end{array}\right)$$
  for some non-zero scalar $\lambda_g$. Therefore, for each $\bar g\in \overline{G}$,
  $X_{\bar g}=(\lambda_g\xi(g))^{-1}) (\varphi_0\xi)*(X_{\bar
  g})$ where $\varphi_0\xi$ is also  an antiautomorphism. In other
  words a fine $\overline{G}$-grading on $I_1$ is respected by
  antiautomorphism $\varphi_0\xi$. Then, by \cite{BShZ},  ${G/H} \cong N_1 \times \ldots
\times N_k $ where $N_i \cong \mathbb{Z}_2 \times \mathbb{Z}_2.$

{\it Case 2.}  Let $\varphi(I_i)=I_i$. Then a fine $\overline{
G}$-grading on each $I_i$ is also compatible with $\varphi$.
Hence, according to \cite{BShZ}, ${G/H} \cong N_1 \times \ldots
\times N_k$ where $N_i \cong \mathbb{Z}_2 \times \mathbb{Z}_2.$

Therefore, $|G|=2\cdot 2^{2l}=2^{2l+1}$, for some natural number
$l$. Moreover, we have that for each $g\in G$, either $g^2=e$ or
$g^4=e$. On the other hand, according to \cite{bsz}, $G=\mathbb{Z}_{n_1}\times \mathbb{Z}_{n_1}\times\ldots\times \mathbb{Z}_{n_k}\times \mathbb{Z}_{n_k}$, $n_i\in \mathbb{N}$. Moreover,
either $n_i=2$ or $n_i=4$. Therefore, $|G|=2^{2r}\cdot
4^{2s}=2^{2r+4s}$, for some natural numbers $r$ and $s$, which is
contradiction.

\end{proof}

In what follows, let  $\tau$ denote an antiautomorphism of
$M_{n,m}(F)$ defined by the formula:
$$ X^{\tau}=\left(\begin{array}{cc}
             A& -B\\
             C&  D
             \end{array}\right)^t,$$
where $A$ and $D$ are matrices of size $n\times n$ and $m\times
m$, respectively, $B$ and $C$ are matrices of size $n\times m$ and
$m\times n$, respectively.

\begin{lemma}
Let $R=M_n=\oplus_{g\in G} R_g$ be a matrix algebra with the
elementary $G$-grading. If $R_e=A_1\oplus A_2$ is the sum of two
simple subalgebras, then there exists $g\in G$, $g\ne e$, such
that $A_1RA_2\subseteq R_g$.
\end{lemma}
\par\medskip
\begin{lemma}
If $R=R_{\bar 0}+R_{\bar 1}$ is a superalgebra with an
antiautomorphism $\varphi$, then for any $x,y\in R$,
$$\varphi(xRy)\subseteq \varphi(y)R\varphi(x).\eqno (2)$$

\end{lemma}

\par\medskip

\begin{lemma}

Let $R=C\otimes D=\oplus_{g\in G} R_g$ be a $G$-graded matrix
superalgebra with an elementary grading on $C$, and a fine grading
on $D$ over an algebraically closed field $F$ of characteristic
not 2. Let $\varphi: R\to R$ be an antiautomorphism on $R$
preserving $G$-grading and $\sigma:R\to R$ be an automorphism of
order 2 of $R$ that defines a superalgebra structure on $R$. Let
also $\varphi$ act as a superinvolution on $R_e$. Then

1) $C_e\otimes I$ is $\varphi$-stable and $\sigma$-stable where
$I$ is the unit of $D$ and hence $\sigma$ induces a $\mathbb{Z}_2$-grading on $C_e$ and $\varphi$ induces a superinvolution $*$
on $C_e$ compatible with $\mathbb{Z}_2$-grading.

2) there are $*$-subsuperalgebras $B_1,\ldots, B_k\subseteq C_e$
such that $C_e=B_1\oplus\ldots\oplus B_k$, and $B_1\otimes
I,\ldots,B_k\otimes I$ are $\varphi$-stable and $\sigma$-stable.

3) $\varphi$ acts on $R_e=C_e\otimes I$ as
$\varphi*X=S^{-1}X^{\tau}S$ where $S=S_1\otimes
I+\ldots+S_k\otimes I$, $S_i\in B_iCB_i$ and
$S_i=\left(\begin{array}{cc}
    I_{s_i}& 0\\
     0 & Q_{r_i}
   \end{array}\right)$ if $B_i$ is of type $M_{s_i,2r_i}(F)$ with
    orthosymplectic superinvolution;
    $S_i=\left(\begin{array}{cc}
            0& I_{s_i}\\
      I_{s_i}&0
          \end{array}\right)
    $ if $B_i$ is of type $M_{s_i,s_i}(F)$ with transpose superinvolution;
    $S_i=\left(\begin{array}{cc}
             0& I_{s_i+r_i}\\
            I_{s_i+r_i}&0
            \end{array}\right)$
if $B_i$ is of type $M_{s_i,r_i}(F)\oplus M^{sop}_{s_i,r_i}(F)$
with exchange superinvolution; $S_i=\left(\begin{array}{cc}
             0& I_{2s_i}\\
            I_{2s_i}&0
            \end{array}\right)$
if $B_i$ is of type $Q(s_i)\oplus Q(s_i)^{sop}$ with exchange
superinvolution.

4) if $e_i$ is the identity of $B_i$, then $D_i=e_i\otimes D$ is
$\varphi$-stable and $\sigma$-stable.

5) the centralizer of $R_e=C_e\otimes I$ in $R$ can be decomposed
as $Z_1D_1\oplus\ldots\oplus Z_kD_k$ where $Z_i=Z'_i\otimes I$,
$Z'_i$ is the center of $B_i$.

\end{lemma}

\begin{proof}

It follows from \cite{bsz} that the identity component $R_e$
equals to $C_e\otimes I$. Since $R_e$  is $\varphi$- and
$\sigma$-stable, both $\varphi$ and $\sigma$ induce a
superinvolution $*$  and a superalgebra structure on $C_e$. Both
structures are compatible with each other.

Since $C_e$ is semisimple, it is the direct sum of simple
subalgebras,
$$ C_e=A_1\oplus\ldots\oplus A_l. $$
If for some $i$, $1\le i\le l$, $\sigma(A_i)=A_j$ where $i\ne j$,
then it is easily seen that $A'_i=A_i+A_j$ is $\sigma$-stable.
Therefore, $C_e$ can be written as a direct sum of $\sigma$-stable
superalgebras,
$$ C_e=A'_1\oplus\ldots\oplus A'_s.$$
Next, if for some $i$, $1\le i\le s$, $(A'_i)^*=A'_j$ where $i\ne
j$, then $B_i=A'_i+A'_j$ is $*$-stable. Finally, $C_e$ can be
written as a direct sum of $*$-simple superalgebras.

$$ C_e=B_1\oplus\ldots\oplus B_k. $$
Now 1), 2) and 3) follows from the classification of involution
simple superalgebras (see Proposition 1).

Next we fix $1\le i\le k$, and consider $R'=(e_i\otimes
I)(C\otimes D)(e_i\otimes I)=e_iCe_i\otimes D$ where $e_i$ is the
identity of $B_i$. Since $\varphi(e_i\otimes I)=e_i\otimes I$ and
$\sigma(e_i\otimes I)=e_i\otimes I$, $R'$ is $\varphi$- and
$\sigma$-stable.

To prove 4), we consider the following three cases.

{\it Case 1.} Let $B_i$ be of the type $M_{r,s}(F)$. Then
$$e_iCe_i=B_i, \eqno (3)$$
and $e_iCe_i\otimes I=B_i\otimes I$. Hence, $e_iCe_i\otimes I$ is
$\varphi$- and $\sigma$-stable. Since $e_i\otimes D$ is a
centralizer of $e_iCe_i\otimes I$, it is also $\varphi$- and
$\sigma$-stable.

{\it Case 2.} Let $B_i=A\oplus A^{sop}$ where $A=M_{r,s}(F)$.
Denote the identity of $A$ by $\varepsilon_i $. Then,
$\varepsilon^*_i$ is the identity of $A^{sop}$, and
$e_i=\varepsilon_i+\varepsilon^*_i$. Notice that $$e_iCe_i\otimes
I=\varepsilon_iC\varepsilon_i\otimes
I+\varepsilon_iC\varepsilon^*_i\otimes
I+\varepsilon^*_iC\varepsilon_i\otimes
I+\varepsilon^*_iC\varepsilon^*_i\otimes I. \eqno (4)$$ Next we
want to prove that both $\varphi$ and $\sigma$ permute the terms
of (4) leaving $e_iCe_i\otimes I$ invariant. Without any loss of
generality we consider just one term of the form
$\varepsilon_iC\varepsilon^*_i\otimes I$. Since
$\varepsilon_iC\varepsilon^*_i\otimes I=(\varepsilon_i\otimes
I)(C\otimes I)(\varepsilon^*\otimes I)$, by (2),
$\varphi(\varepsilon_iC\varepsilon^*_i\otimes
I)\subseteq(\varepsilon_i\otimes I)(C\otimes
D)(\varepsilon^*\otimes I)=\varepsilon_iC\varepsilon^*_i\otimes D$
and $\sigma(\varepsilon_iC\varepsilon^*_i\otimes
I)\subseteq(\varepsilon_i\otimes I)(C\otimes
D)(\varepsilon^*\otimes I)=\varepsilon_iC\varepsilon^*_i\otimes
D$.

By Lemma 2, there exists $g\in G$, $g\ne e$ such that
$\varepsilon_i C\varepsilon^*_i\subseteq C_g$. Hence,
$\varepsilon_i C\varepsilon^*_i\otimes I\subseteq  R_g$.
Consequently, $\varphi(\varepsilon_i C\varepsilon^*_i\otimes
I)\subseteq  R_g$ and $\sigma(\varepsilon_i
C\varepsilon^*_i\otimes I)\subseteq  R_g.$ Next we take a
homogeneous $x\in \varepsilon_i C\varepsilon^*_i$ of degree $g$
and a homogeneous $y\in D$ of degree $h$ such that $x\otimes y \in
R_g$. Then $\mbox{deg}\, (x\otimes y)=gh=g$, $h=e$. This implies
$y=\lambda I$, $\lambda\in F$ for any $x\otimes y\in R_g\cap
\varepsilon_i C\varepsilon^*_i\otimes D$.  It follows that $
R_g\cap \varepsilon_i C\varepsilon^*_i\otimes D\subseteq
\varepsilon_i C\varepsilon^*_i\otimes I$.

As a result,
 $\varphi(e_iCe_i\otimes I)=e_iCe_i\otimes I $ and $\sigma(e_iCe_i\otimes I)=e_iCe_i\otimes
 I$, that is, $e_iCe_i\otimes I$ is
$\varphi$- and $\sigma$-stable. From the decomposition
$R'=e_iCe_i\otimes D$ it follows that $e_i\otimes D$, the
centralizer of $e_iCe_i\otimes I$ in $R'$, is $\varphi$- and
$\sigma$-stable.

{\it Case 3.} Let  $B_i=Q(s_i)\oplus Q(s_i)^{sop}$. Since
$Q(s_i)=I_1\oplus I_2$ where $I_1$, $I_2$ are simple ideals
isomorphic to  $M_{s_i}(F)$, $B_i=(I_1\oplus I_2)\oplus
(I^*_1\oplus I^*_2)$. Let $\varepsilon_i$, ${\hat \varepsilon}_i$,
$\varepsilon^*_i$, $\hat \varepsilon^*_i$  be the identities of
$I_1,$ $I_2$, $I^*_1$, $I^*_2$, respectively. Then we notice that
$\sigma(\varepsilon_i\otimes I)={\hat \varepsilon}_i\otimes I$. We
have that
$$e_i=\varepsilon_i+\varepsilon^*_i+{\hat \varepsilon}_i+{\hat
\varepsilon}^*_i. \eqno (5)$$

Therefore, $e_iCe_i\otimes I=N_1\otimes I +N_2\otimes I+N_3\otimes
I+N_4\otimes I$ where
$N_1=\varepsilon_iC\varepsilon_i+\varepsilon^*_iC\varepsilon_i+\varepsilon_iC\varepsilon^*_i+\varepsilon^*_iC\varepsilon^*_i$,
$N_2=\varepsilon_iC\hat\varepsilon_i+\varepsilon^*_iC\hat\varepsilon_i+\varepsilon_iC\hat\varepsilon^*_i+\varepsilon^*_iC\hat\varepsilon^*_i$,
$N_3=\hat\varepsilon_iC\varepsilon_i+\hat\varepsilon^*_iC\varepsilon_i+\hat\varepsilon_iC\varepsilon^*_i+\hat\varepsilon^*_iC\varepsilon^*_i$,
and
$N_4=\hat\varepsilon_iC\hat\varepsilon_i+\hat\varepsilon^*_iC\hat\varepsilon_i+\hat\varepsilon_iC\hat\varepsilon^*_i+\hat\varepsilon^*_iC\hat\varepsilon^*_i$.

Arguing in the same way as in the second case, we can prove that
$\varphi(N_1\otimes I)=N_1\otimes I$ and $\varphi(N_4\otimes
I)=N_4\otimes I$. Now we consider $N_2$ and $N_3$. Suppose that
the elementary grading on $B_iCB_i$ induced from $C$ is defined by
$(g_1,g_2,g_3,g_4)$. It is easy to see that $\deg
(\varepsilon_iC\hat\varepsilon_i)=g^{-1}_1g_3,$ $\deg
(\varepsilon_iC\hat\varepsilon^*_i)=g^{-1}_1g_4,$ $\deg
(\varepsilon^*_iC\hat\varepsilon_i)=g^{-1}_2g_3,$ $\deg
(\varepsilon^*_iC\hat\varepsilon^*_i)=g^{-1}_2g_4,$ $\deg
(\hat\varepsilon_iC\varepsilon_i)=g^{-1}_3g_1,$ $\deg
(\hat\varepsilon^*_iC\varepsilon_i)=g^{-1}_4g_1,$ $\deg
(\hat\varepsilon^*_iC\varepsilon^*_i)=g^{-1}_4g_2,$ $\deg
(\hat\varepsilon_iC\varepsilon^*_i)=g^{-1}_3g_2.$

Let us take, for example, the first term
$\varepsilon_iC\hat\varepsilon_i$ of $N_2$. Then
$\varphi(\varepsilon_iC\hat\varepsilon_i\otimes I)\subseteq
R_{g^{-1}_1g_3}$. On the other hand, by (2),
$\varphi(\varepsilon_iC\hat\varepsilon_i\otimes I)\subseteq
\hat\varepsilon^*_iC\varepsilon^*_i\otimes D\subseteq
C_{g^{-1}_4g_2}\otimes D$. If we take $x\in C_{g^{-1}_4g_2}$ and a
homogeneous $y\in D$ of degree $h$ such that $\deg (x\otimes
y)=g^{-1}_1g_3$, then $g^{-1}_1g_3=g^{-1}_4g_2h$, that is,
$h=g^{-1}_1g_4g^{-1}_2g_3$. This implies that
$$\varphi(\varepsilon_iC\hat\varepsilon_i\otimes I)\subseteq
\hat\varepsilon^*_iC\varepsilon^*_i\otimes D_h. \eqno(6)$$
Similarly, we can show that for each term of $N_2$ there should be
$D_h$ on the right-hand side of (6). Hence $\varphi(N_2\otimes
I)=N_3\otimes D_{h}$.

Next we take the first term $\hat\varepsilon_iC\varepsilon_i$ of
$N_3$. Likewise we can show that $\varphi(\hat
\varepsilon_iC\varepsilon_i\otimes I)\subseteq
\varepsilon^*_iC\hat\varepsilon^*_i\otimes D_{h^{-1}}$. Hence for
each term of $N_3$ we have $D_{h^{-1}}$ on the right hand side,
and $\varphi(N_3\otimes I)=N_2\otimes D_{h^{-1}}$.

Notice that the centralizer of $(N_2+N_3)\otimes I$ in $R'$ is
$e_i\otimes D$. Hence, the centralizer of
$\varphi((N_2+N_3)\otimes I)=N_3\otimes D_{h^{-1}}+N_2\otimes D_h$
in $R'$ is $\varphi(e_i\otimes D)$. Next we take any $x\otimes y$
in $\varphi(e_i\otimes D)$. Since $x\otimes y$ lies in the
centralizer of $N_3\otimes D_{h^{-1}}+N_2\otimes D_h$, $x\otimes
y$ commutes with each element in $N_3\otimes D_{h^{-1}}$ and
$N_2\otimes D_h$. Therefore, $x$ commutes with any matrix in $N_3$
and $N_2$. Direct computations show that $x=e_i$, and
$\varphi(e_i\otimes D)=e_i\otimes K$ where $K$ is a subspace of
$D$. By dimension arguments, $K=D$, and $\varphi(e_i\otimes D)=
e_i\otimes D$.

To prove that $e_i\otimes D$ is $\sigma$-stable, we represent
$e_iCe_i\otimes I$ as follows: $e_iCe_i\otimes I=N'_1\otimes I
+N'_2\otimes I+N'_3\otimes I+N'_4\otimes I$ where
$N'_1=\varepsilon_iC\varepsilon_i+\hat\varepsilon_iC\varepsilon_i+\varepsilon_iC\hat\varepsilon_i+\hat\varepsilon_iC\hat\varepsilon_i$,
$N'_2=\varepsilon_iC\varepsilon^*_i+\hat\varepsilon_iC\varepsilon^*_i+\varepsilon_iC\hat\varepsilon^*_i+\hat\varepsilon_iC\hat\varepsilon^*_i$,
$N'_3=\varepsilon^*_iC\varepsilon_i+\hat\varepsilon^*_iC\varepsilon_i+\varepsilon^*_iC\hat\varepsilon_i+\hat\varepsilon^*_iC\hat\varepsilon_i$,
and
$N'_4=\varepsilon^*_iC\varepsilon^*_i+\hat\varepsilon^*_iC\varepsilon^*_i+\varepsilon^*_iC\hat\varepsilon^*_i+\hat\varepsilon^*_iC\hat\varepsilon^*_i$.
In the same way as above, we can show that $e_i\otimes I$ is
$\sigma$-stable.

 Hence 4) is proved.

To prove 5) we note that the centralizer $Z$ of $C_e$ in $C$ is
equal to $Z'_1\oplus\ldots\oplus Z'_k$ where $Z'_i$ is the center
of $B_i$ and the centralizer of $R_e$ in $R$ coincides with
$Z\otimes D=Z_1D_1\oplus\ldots\oplus Z_kD_k$ where
$Z_i=Z'_i\otimes I$ and $D_i=e_i\otimes D$.

Our proof is complete.

\end{proof}

\begin{proposition}
Let $G$ be a finite abelian group, and $R$  a superalgebra of type
$M_{n,m}(F)$ over an algebraically closed field $F$ of
characteristic not 2. Suppose that $\varphi$ is an
antiautomorphism of $R$ that preserves a $G$-grading of $R$. Then
there exists an automorphism $\psi$ of $R$ preserving the
$G$-grading of $R$ such that $\psi$ commutes with $\varphi$ and
$\psi^4=\varphi^4$.

\end{proposition}

\begin{proof}
It is easy to check that the $\varphi$-action on $R$ is defined by
$$
\varphi * X=\Phi^{-1}X^{\tau}\Phi
$$
for some matrix $\Phi$, and  $X^{\tau}=\left(\begin{array}{cc}
             A& -B\\
             C&  D
             \end{array}\right)^t$. First let $X\in R_e$. Consider the
decomposition $C_e=B_1\oplus\ldots\oplus B_k$ found in the
previous lemma. Then $X=X_1\otimes I+\ldots+ X_k\otimes I$ with
$X_i\in B_i, 1\le i \le k$. Then $\varphi$ acts on $X$ as
$$
\varphi * X=S^{-1}X^{\tau}S,
$$
where $S$ as in 3) of Lemma 4. Hence the matrix $\Phi S^{-1}$
commutes with $X^{\tau}$ for any $X\in R_e$, that is $\Phi S^{-1}$
is an element of the centralizer of $R_e$ in $R$. Hence,  we
obtain
$$
\Phi= S_1Y_1\otimes Q_1+\ldots+ S_kY_k\otimes Q_k \eqno (7)
$$
where $Q_i\in D,Y_i\in Z_i', 1\le i\le k$. Compute now the action
of $\varphi^4$ on an arbitrary $X\in R$:
$$
\varphi^4*X=((\Phi^{-1})^t\Phi)^{-1})^2X((\Phi^{-1})^t\Phi)^2
$$

Set $P=((\Phi^{-1})^t\Phi)^2$. We need to  show that there exists
an inner automorphism $\psi$ such that $\psi^4* X=P^{-1}XP$ for
all $X\in R$. Note that for any $T_i,T_i'\in B_iCB_i$ and
$Q_i,Q_i'\in D$, $i=1,\ldots,k$, the relation
$$
\left(\sum_i T_i\otimes Q_i\right)\left(\sum_i T_i'\otimes
Q_i'\right)= \sum_i T_iT_i'\otimes Q_iQ_i'
$$
holds.

We compute the value of $P$:
$$P=((\Phi^{-1})^t\Phi )^2=\sum_{i=1}^k
({(Y^t_iS^t_i)}^{-1}S_iY_i)^2\otimes ({(Q^t_i)}^{-1}Q_i)^2= $$
$$\sum_i ({(S^t_i)}^{-1}{(Y^t_i)}^{-1}S_iY_i)^2\otimes
\,({(Q^t_i)}^{-1}Q_i)^2.\eqno (8)
$$

\begin{lemma}
All $Q_i$  satisfy ${(Q^t_i)}^{-1}Q_i=\pm I$.
\end{lemma}
Obviously it is sufficient to prove the relation
$$ e_i\otimes (Q^t_i)^{-1}Q_i=\pm e_i\otimes I$$
in $D_i=e_i\otimes D$. Recall that $D_i$ is $\varphi$- and
$\sigma$-stable. Moreover, $D_i$ is $G$-graded algebra with a fine
$G$-grading compatible with $\varphi$ and $\sigma$. Therefore,
this is $G$-graded superalgebra with a fine $G$-grading respected
by $\varphi$. According to Theorem 3, $D_i$ cannot be non-trivial.
Therefore, $D_i$ is a trivial superalgebra, that is, $D_i\subseteq
R_{\bar 0}$, and $\tau$ acts on $D_i$ as a usual transpose. For
any $X\in D$ we have
$$\varphi*(e_i\otimes X)=\Phi^{-1}(e_i\otimes
X)^{\tau}\Phi=(S_iY_i)^{-1}(e_i)(S_iY_i)\otimes
(Q^{-1}_iX^tQ_i)=e_i\otimes Q^{-1}_iX^tQ_i$$ i.e. $\varphi$-action
induces an antiautomorphism $e_i\otimes X\to e_i\otimes
Q^{-1}_iX^tQ_i$ on $D_i$. Arguing in the same way as in Lemma 6.5
(see \cite{bs}) we can conclude that $e_i\otimes
(Q^t_i)^{-1}Q_i=\pm e_i\otimes I$
\end{proof}
Now we compute $({(Y^t_iS^t_i)}^{-1}S_iY_i)^2$. If $B_i$ is simple
then $Y_i$ is a scalar matrix and $((S^t_i)^{-1}S_i)^2=I$ by Lemma
4. If $B_i$ is of type $M_{s_i,r_i}(F)\oplus
M^{sop}_{s_i,r_i}(F)$, then
$$Y_i=\left(\begin{array}{cc}
          \lambda I& 0\\
                  0& \mu I
       \end{array}\right),\quad
  S_i=\left(\begin{array}{cc}
                  0& I\\
                  I &0
            \end{array}\right)
$$

and $({(Y^t_iS^t_i)}^{-1}S_iY_i)^2=\left(\begin{array}{cc}
          (\frac{\lambda}{\mu})^2 I& 0\\
                  0& (\frac{\mu}{\lambda})^2 I
       \end{array}\right)=
\left(\begin{array}{cc}
          (\gamma)^2 I& 0\\
                  0& (\gamma^{-1})^2 I
       \end{array}\right)$ where $\gamma=\frac{\lambda}{\mu}$.
If $B_i$ if of type $Q(s_i)\oplus Q(s_i)^{sop}$, then

$$Y_i=\left(\begin{array}{cccc}
        \alpha I& 0& 0& 0\\
        0& \alpha_1 I& 0& 0\\
        0&      0& \beta I & 0\\
        0&      0&   0  & \beta_1 I
        \end{array}\right),\quad
       S_i=\left(\begin{array}{cccc}
        0& 0& I& 0\\
        0& 0& 0& I\\
        I&  0& 0& 0\\
        0&  I&   0  & 0
        \end{array}\right),
$$
and
$$({(Y^t_iS^t_i)}^{-1}S_iY_i)^2=
\left(\begin{array}{cccc}
        (\beta^{-1}\alpha)^2 I& 0& 0& 0\\
        0& (\alpha_1\beta^{-1}_1)^2 I& 0& 0\\
        0&      0& (\alpha^{-1}\beta)^2 I & 0\\
        0&      0&   0  & (\alpha^{-1}_1\beta_1)^2 I
        \end{array}\right)=
$$
$$\left(\begin{array}{cccc}
        \gamma^2 I& 0& 0& 0\\
        0& \mu^2 I& 0& 0\\
        0&      0& (\gamma^{-1})^2 I & 0\\
        0&      0&   0  & (\mu^{-1})^2 I
        \end{array}\right)
$$
where $\gamma=\beta^{-1}\alpha$, and $\mu=\alpha_1\beta^{-1}_1$.

We have proved that $P=(P_1+\ldots+P_k)\otimes I$ where $P_1\in
B_1,\ldots,P_k\in B_k$ and $P_i$ has one of the forms: $P_i=I$,
$P_i=\left(\begin{array}{cc}
          (\gamma)^2 I& 0\\
                  0& (\gamma^{-1})^2 I
       \end{array}\right)$, and
       $P_i=\left(\begin{array}{cccc}
        \gamma^2 I& 0& 0& 0\\
        0& \mu^2 I& 0& 0\\
        0&      0& (\gamma^{-1})^2 I & 0\\
        0&      0&   0  & (\mu^{-1})^2 I
        \end{array}\right).$ Now we present a matrix
$T=(T_1+\ldots+T_k)\otimes I$, $T_1\in B_1,\ldots, T_k\in B_k$
such that $T^4_i=P_i$ for all $i$ and hence $T^4=P$. In case
$P_i=I$ we take $T_i=I$. If $P_i=\left(\begin{array}{cc}
          (\gamma)^2 I& 0\\
                  0& (\gamma^{-1})^2 I
       \end{array}\right)$, then we take
$T_i=\left(\begin{array}{cc}
          \gamma_1 I& 0\\
                  0& \gamma^{-1}_1 I
       \end{array}\right)$ where $\gamma^2=\gamma^4_1$. If
   $P_i=\left(\begin{array}{cccc}
        \gamma^2 I& 0& 0& 0\\
        0& \mu^2 I& 0& 0\\
        0&      0& (\gamma^{-1})^2 I & 0\\
        0&      0&   0  & (\mu^{-1})^2 I
        \end{array}\right)$, we take
 $T_i=\left(\begin{array}{cccc}
        \gamma_1 I& 0& 0& 0\\
        0& \mu_1 I& 0& 0\\
        0&      0& \gamma^{-1}_1 I & 0\\
        0&      0&   0  & \mu^{-1}_1 I
        \end{array}\right)$
 where $\gamma^2=\gamma^4_1$, and $\mu^2=\mu^4_1$.
Note that $T\in R_e$, hence the map $\psi:X\to T^{-1}XT$ is an
inner automorphism preserving $G$-grading. Moreover, since
$T^4=P$, $\psi^4=\varphi^4$.

Now we need to check that $\psi$ and $\varphi$ commute. Direct
computations show that $\varphi\psi=\psi\varphi$ if and only if
$$ T^{\tau}\Phi T=\lambda \Phi,\quad (\Phi^{-1}T^{\tau}\Phi)T=\lambda I, \eqno (9) $$
for some scalar $\lambda$. Since $T=T_1\otimes I+\ldots+T_k\otimes
I$ where $T_i\in B_i$, $\Phi^{-1}T^{\tau}\Phi=\varphi*T=
S^{-1}T^{\tau}S$ (see Lemma 4). If $B_i$ is simple, then $T_i=I$
and $S^{-1}_iT_i^{\tau}S_i=T_i=I$. If $B_i=A\oplus A^{sop}$, then
the restriction of $\varphi$ to $B_i$ acts as the exchange
superinvolution, and
$$ S^{-1}_iT_i^{\tau} S_i=\left(\begin{array}{cccc}
        \gamma^{-1}_1 I& 0& 0& 0\\
        0& \mu^{-1}_1 I& 0& 0\\
        0&      0& \gamma_1 I & 0\\
        0&      0&   0  & \mu_1 I
        \end{array}\right)
$$
for $T_i=\left(\begin{array}{cccc}
        \gamma_1 I& 0& 0& 0\\
        0& \mu_1 I& 0& 0\\
        0&      0& \gamma^{-1}_1 I & 0\\
        0&      0&   0  & \mu^{-1}_1 I
        \end{array}\right)$, or
$$ S^{-1}_iT_i^{\tau} S_i=\left(\begin{array}{cc}
          \gamma^{-1}_1 I& 0\\
                  0& \gamma_1 I
       \end{array}\right)$$
for $T_i=\left(\begin{array}{cc}
          \gamma_1 I& 0\\
                  0& \gamma^{-1}_1 I
       \end{array}\right)$.
In both cases (9) holds with  $\lambda=1$ and thus the proof is
complete.

\par\medskip

\section{Main results}

In this section we describe group gradings compatible with
superinvolution of involution simple superalgebras which are not
simple as superalgebras. Notice that these results depend on the
classification of  gradings by a finite abelian group on matrix
algebras \cite{bsz}, involution gradings on matrix algebras
\cite{BShZ}, \cite{ant}, involution gradings on involution simple
algebras \cite{bg}, and group gradings on simple superalgebras
\cite{BSh}. Finally, the superinvolution gradings on $M_{n,m}(F)$
have been described in \cite{TT}.

We start the following general result.
\begin{lemma}\label{l5}
Let $R$ be a simple superalgebra of Example 3, that is, $R=A\oplus A^{sop}$ where $A$ is a simple
superalgebra, and $*$ denote the ordinary exchange involution.
If $\varphi$ is an automorphism of $R$ that commutes with $*$,
then there exists a linear mapping $\varphi_0:A\to A$ such that
one of the following cases holds:
\begin{enumerate}
\item[Type 1:] $\varphi((x,y))=(\varphi_0(x),\varphi_0(y))$, and $\varphi_0$ is an automorphism
of $A$.
\item[Type 2:] $\varphi((x,y))=(\varphi_0(y),\varphi_0(x))$, and $\varphi_0$ is an antiautomorphism
of $A$.
\end{enumerate}
\end{lemma}
\begin{proof}
Since $R=A\oplus A^{sop}$, we will represent an arbitrary element
of a superalgebra $R$ as a pair of elements from $A$, i.e. $(x,y)$
where $x,y\in A$.  We also recall that $A=A_{\bz}+A_{\bo}$. If $\varphi$
is an automorphism of $R$ that commutes with $*$, then the
following two cases may occur:

1. $\varphi(A)=A$, $\varphi(A^{sop})= A^{sop}$.  Then, there
exist two linear mappings $\varphi_0,\varphi_1: A\to A $ such that
$\varphi((x,y))=(\varphi_0(x),\varphi_1(y))$. Now $\varphi$
commutes with the involution $*$. Hence
$$(\varphi_1(y),\varphi_0(x))=(\varphi_0(x),\varphi_1(y))^*=(\varphi((x,y))^*$$
$$=\varphi((x,y))^*=\varphi((y,x))=(\varphi_0(y),\varphi_1(x)).$$
Hence, $\varphi_0=\varphi_1$. Thus $\varphi$ is completely defined
by $\varphi_0: A\to A$,
$\varphi((x,y))=(\varphi_0(x),\varphi_0(y))$. Next, for any
homogeneous $x,y\in A$,
\begin{eqnarray*}
\varphi((xy,0))&=&(\varphi((0,xy)))^*=((-1)^{|x||y|}\varphi((0,y))\varphi((0,x)))^*\\
&=&(-1)^{|x||y|}((0,\varphi_0(y))\cdot(0,\varphi_0(x)))^*\\&=&(-1)^{|x||y|}(-1)^{|\varphi_0(x)||\varphi_0(y)|}
 (0,\varphi_0(x)\varphi_0(y))^*\\&=&(\varphi_0(x)\varphi_0(y),0).\end{eqnarray*}
\noindent Hence $\varphi_0$ is indeed an automorphism of
$A$, and we have a Type I grading.

2. $\varphi(A)=A^{sop}$, $\varphi(A^{sop})=A$. Again there exist
two linear mappings $\varphi_0,\varphi_1: A\to A$ such that $\varphi((x,y))=
(\varphi_0(y),\varphi_1(x))$. Since $\varphi$  commutes with
the involution, we must have
$$(\varphi_1(x),\varphi_0(y))=\varphi((x,y))^*=\varphi((y,x))=
(\varphi_0(x),\varphi_1(y)).$$ Again, as before
$\varphi_0=\varphi_1$. Now let $x,y$ be homogeneous elements from
$A$. Therefore,
$(\varphi_0(xy),0)=\varphi((0,xy))=\varphi(((-1)^{|x||y|}(0,y)\cdot
(0,x)))=(-1)^{|x||y|}\varphi((0,y))\cdot\varphi((0,x))=
(-1)^{|x||y|}(\varphi_0(y),0)\cdot (\varphi_0(x),0)$. It follows
that $$\varphi_0(xy) = (-1)^{|x||y|}\varphi_0(y)\varphi_0(x),$$ and we have a Type II grading.

\end{proof}

\par\medskip

 By a \emph{Type I} involution grading of a superalgebra
$R_g=A\oplus A^{sop}$, as above, we understand a grading in which $A$ is
a graded subspace, that is, $A=\bigoplus_{g\in G}(R_g\cap A)$. In
this case also $$A^{sop} = A^\ast = \bigoplus_{g\in
G}(R_g^{\ast}\cap A^{\ast})=\bigoplus_{g\in G}(R_g\cap A^{sop}),$$
so that $A^{sop}$ is also graded. Then there is a $G$-grading on
$A$, hence on $A^{sop}$, as in Remark \ref{r1}, such that
$R_g=A_g\oplus A_g^{sop}$.

\begin{theorem}\label{tm1} Let $G$ be a finite abelian group and $F$
an algebraically closed field of characteristic 0 or coprime to
the order of $G$. Then any $G$-grading of $R=A\oplus A^{sop}$
where $A=M_{n,m}(F)$, with the ordinary exchange superinvolution
$*$ compatible with $G$-grading has one of the following types:

\emph{Type I:} $R_g=A_g\oplus A_g^{sop}$, for  a $G$-grading of
$A=\oplus_{g\in G} A_g,$

\emph{Type II:} $R_g=\{(x,x^{\dagger})|\,x\in A_g\}\oplus
\{(x,-x^{\dagger})|\,x\in A_{gh}\}$, for a $\dagger$-involution
grading $A=\oplus_{g\in G} A_g$ where $\dagger$ is a graded
superinvolution on $A$, $h\in G$, $o(h)=2$.

\emph{Type III:} $R_g=\{(x,x^{\dagger})\,|\, x\in A_g\cap
A_+\}\oplus \{(x,-ix^{\dagger})\,|\, x\in A_{gh}\cap A_-\}\oplus
\{(x,-x^{\dagger})\,|\, x\in A_{gh^2}\cap A_+\}\oplus
\{(x,ix^{\dagger})\,|\, x\in A_{gh^3}\cap A_-\}$  where $h$ is an
element of order 4 in $G$, $\dagger$ is an antiautomorphism of
order 4 on $A$, $A=\bigoplus_{g\in G}A_g$ is a $\dagger$-grading
on $A$, $A_{+}$, $A_{-}$ are symmetric and skew-symmetric elements
of $A$ with respect to $\dagger^2$.

\end{theorem}

\begin{proof}
If $\wg$ acts on $R$ by the automorphisms of Type 1 only, we
arrive at Type I gradings described just before the statement of
this theorem. Now let $\wg$ act on $R$ by the automorphisms of
both Type 1 and Type 2, and  $\alpha:\wg\to \mbox{Aut}\,R$ the
homomorphism accompanying our grading. Let $\Lambda$ stand for the
set of all $\chi\in\wg$ that act on $R$ by the automorphisms of
Type 1. As earlier, $\Lambda$ is a subgroup of index 2 in $\wg$.
Choose $\xi\in\wg$, such that $\alpha(\xi)=\varphi$ is an automorphism of
Type 2, $\wg=\Lambda\cup\Lambda \xi$.

Next we assume that there exists an automorphism $\psi$ of  Type 1
such that $\psi^2=\varphi^2$, and $\psi$ commutes with
$\alpha(\widehat G)$. Then we can apply the Exchange
Theorem. For this, we consider two gradings of $R$. The first is
our original one defined by $\alpha$. The second one is defined by
a new homomorphism $\beta$ such that
$\beta|_{\Lambda}=\alpha|_{\Lambda}$, $\beta(\xi)=\psi$. It is
easily seen that $\beta$ is indeed a homomorphism. Now by the
Exchange Theorem there exists a grading by a subgroup
$H=\Lambda^{\bot}=\{e,h\}$, corresponding to the action of
$\gamma=\alpha\beta^{-1}$. Moreover,
$\gamma(\wg)=\{id,\varphi\psi^{-1}\}$. Denote
$\omega=\varphi\psi^{-1}$. Then,
$\omega((x,y))=(\omega_0(y),\omega_0(x))$, $\omega_0^2=id$ and
$\omega_0(xy)=(-1)^{|x||y|}\omega_0(y)\omega_0(x)$. Therefore,
$\omega_0$ is an involution on $M_{n,m}(F)$ which we denote by
$\omega_0(x)=x^{\dagger}$. Since the grading defined by $\beta$ is
a grading of the Type I, it follows from the first part of the
proof of this theorem that $\overline{R}_g=A_g\oplus A_g^{sop}$
where $A=\oplus_{g\in G}A_g$ is a $\dagger$-grading of $R$. By the
Exchange Theorem there exists an element $h$ of $G$ of order 2
such that  $R_g=\overline{R}_g\cap R^{(e)}+\overline{R}_{gh}\cap
R^{(h)}$. Here,
$$ R^{(e)}=\{(x,y)|\omega((x,y))=(x,y)\}=$$
$$ \{(x,y)|(\omega_0(y),\omega_0(x))=(x,y)\}=
\{(x,x^{\dagger})|\,x\in A\}.$$

\noindent Also

$$ R^{(h)}=\{(x,y)|\omega((x,y))=-(x,y)\}=\{(x,-x^{\dagger})|x\in A\}.$$ This allows us to
write $R_g=\{(x,x^{\dagger})|\,x\in
A_g\}\cup\{(x,-x^{\dagger})|\,x\in A_{gh}\}$.

Now we consider the remaining case when there is no automorphism
$\psi$ of $R$ of Type 1 such that $\psi^2=\varphi^2$ and $\psi$
commutes with $\alpha({\widehat G})$. Let $\Lambda_1$ denote the
set of all $\eta\in\widehat G$ for which there exists an
automorphism $\tau$ of Type 1 such that $\alpha(\eta)=\tau^2$ and
$\tau$ commutes with $\alpha(\widehat G)$. Clearly, $\Lambda_1$ is
a subgroup of $\Lambda$. Moreover, since $\eta^2\in \Lambda_1$ for
each $\eta\in\Lambda$, this subgroup has index 2 in $\Lambda$, and
therefore, has index 4 in $\widehat G$. By our assumption
$\xi^2\notin \Lambda_1$. Hence $\widehat
G=\Lambda_1\cup\Lambda_1\xi\cup\Lambda_1\xi^2\cup\Lambda_1\xi^3$.
Next we can write  $\varphi((x,y))=(\varphi_0(y),\varphi_0(x))$
where $\varphi_0$ is an antiautomorphism of $A$ that commutes with
$\alpha(\Lambda_1)$. By Proposition 3, there exists an
automorphism $\psi_0$ of $A$ such that $\psi_0^4=\varphi_0^4$,
$\psi_0\varphi_0=\varphi_0\psi_0$, and $\psi_0$ commutes with
$\alpha(\Lambda_1)$.
 Set
$\psi((x,y))=(\psi_0(x),\psi_0(y))$. Obviously,
$\psi^4((x,y))=(\psi_0^4(x),\psi_0^4(y))=(\varphi_0^4(x),\varphi_0^4(y))=
\varphi^4((x,y))$. Besides,
$$\psi\varphi((x,y))=\psi((\varphi_0(y),\varphi_0(x)))=(\psi_0\varphi_0(y),\psi_0\varphi_0(x))=$$
$$(\varphi_0\psi_0(y),\varphi_0\psi_0(x))=\varphi((\psi_0(x),\psi_0(y)))=\varphi\psi((x,y)).$$
This implies that $\psi\varphi=\varphi\psi$. Moreover, $\psi$
commutes with $\alpha(\Lambda_1)$.  Next we consider a new
homomorphism $\beta:{\widehat G}\to\mbox{Aut}\,R$ defined as
follows: $\beta(\xi^k\eta)=\psi^k\alpha(\eta)$ for $k=0,1,2,3$.
Let also $R=\bigoplus_{g\in G}\overline{R}_g$ be a $G$-grading
defined by $\beta$. Since this is a Type I grading, ${\overline
R}_g=A_g\oplus A^{sop}_g$ for some $G$-grading of
$A=\oplus_g\,A_g$.

 This allows us to apply Exchange Theorem, in which
$\Lambda^{\perp}=H=\{e,h,h^2,h^3\}$. The homomorphism
$\gamma:\wg\rightarrow\aut{R}$ defined by
$\gamma(\chi)=\alpha^{-1}(\chi)\beta(\chi)$ defines a grading
$R=R^{(e)}\oplus R^{(h)}\oplus R^{(h^2)}\oplus R^{(h^3)}$ and
$$
R_g=\overline{R}_g\cap R^{(e)}\oplus \overline{R}_{gh}\cap
R^{(h)}\oplus \overline{R}_{gh^2}\cap R^{(h^2)}\oplus
\overline{R}_{gh^3}\cap R^{(h^3)}.
$$

Let us consider $\theta=\gamma(\xi)$. Then the respective
$\theta_0$ is an antiautomorphism of $A$ of order 4 which we
denote by $\dagger$. Notice that $\theta^2$ is an automorphism of
order 2.  Let $A_{+}=\{x\in A|\,\, \theta^2_0(x)=x\}$ and
$A_{-}=\{x\in A|\,\, \theta^2_0(x)=-x\}$. Direct computations show
that  $\overline{R}_g\cap R^{(e)}=\{(x,x^{\dagger})\,|\, x\in
A_g\cap A_+\}$,  $\overline{R}_{gh}\cap
R^{(h)}=\{(x,-ix^{\dagger})\,|\, x\in A_{gh}\cap A_-\}$,
$\overline{R}_{gh^2}\cap R^{(h^2)}=\{(x,-x^{\dagger})\,|\, x\in
A_{gh^2}\cap A_+\}$ and $\overline{R}_{gh^3}\cap
R^{(h^3)}=\{(x,ix^{\dagger})\,|\, x\in A_{gh^3}\cap A_-\}$.

 The proof is now complete.

\end{proof}

\noindent{\bf Example}. Let us consider $\mathbb{Z}_4=\{\pm 1,\pm
i\}$-grading of $R=A\oplus A^{sop}$, $A\cong M_{n,m}(F)$ induced
by $\varphi*(X,Y)=(Y^{\tau},X^{\tau})$. Clearly, $\varphi$ is of
order 4. Direct computations show that

$$ R_1=\left\{ \left(
         \left[\begin{array}{cc}
                  A& 0\\
                  0& D
                  \end{array}
              \right],\,\,
          \left[\begin{array}{cc}
                  A^t& 0\\
                  0& D^t
                  \end{array}
              \right]
             \right)\right\},\quad
R_{-1}=\left\{ \left(
         \left[\begin{array}{cc}
                  A& 0\\
                  0& D
                  \end{array}
              \right],\,\,
          \left[\begin{array}{cc}
                  -A^t& 0\\
                  0& -D^t
                  \end{array}
              \right]
             \right)\right\},
$$

$$ R_i=\left\{ \left(
         \left[\begin{array}{cc}
                  0& B\\
                  C& 0
                  \end{array}
              \right],\,
          \left[\begin{array}{cc}
                  0& -iC^t\\
                  iB^t& 0
                  \end{array}
              \right]
             \right)\right\},\,
R_{-i}=\left\{ \left(
         \left[\begin{array}{cc}
                  0& B\\
                  C& 0
                  \end{array}
              \right],\,
          \left[\begin{array}{cc}
                  0& iC^t\\
                  -iB^t& 0
                  \end{array}
              \right]
             \right)\right\},
$$
\noindent where $A,\,B,\,C,\,D$ are any matrices of appropriate
orders. This is in fact a grading of Type III for $A=A_e$ (a
trivial grading) and $h=-i$.

\par\medskip
\par\medskip

\begin{lemma}\label{ltm21}
Let $A=B+tB$ where $B\cong M_n(F)$ be an associative superalgebra
of type $Q(n)$, and $\psi$ an automorphism of $A$. Then there
exists an automorphism $\psi_0$ of $B$ such that for any $X+tY\in
R$, either $\psi(X+tY)=\psi_0(X)+t\psi_0(Y)$ or
$\psi(X+tY)=\psi_0(X)-t\psi_0(Y)$.
\end{lemma}

\begin{proof}
Let us consider $A=B+tB$ with a $\mathbb{Z}_2$-grading $A_{\bar
0}=B$ and $A_{\bar 1}=tB$. Then both $B$ and $tB$ are invariant
subspaces under the action of $\psi$. Namely, there exists two
linear mappings $\psi_0,\psi_1: B\to B$ such that for any $X+tY\in
R$, $\psi(X+tY)=\psi_0(X)+t\psi_1(Y)$. If we use that $\psi$ is an
automorphism, we can easily derive the following relations:

$$\psi_0(X_1X_2)=\psi_0(X_1)\psi_0(X_2) \eqno(10)$$
$$ \psi_1(XY)=\psi_0(X)\psi_1(Y)\eqno(11)$$
 $$\psi_1(YX)=\psi_1(Y)\psi_0(X)\eqno(12)$$
 $$\psi_0(Y_1Y_2)=\psi_1(Y_1)\psi_1(Y_2).\eqno(13)$$
where all $X_1,X_2,X, Y_1,Y_2,Y\in B$. It follows from (10) that
$\psi_0$ is an automorphism of $B$. Now in (11) and (12) we set
$Y=I$, the identity matrix, then we obtain
$\psi_0(X)\psi_1(I)=\psi_1(I)\psi_0(X)$ for all $X\in B$. It
follows then that $\psi_1(I)$ is a scalar matrix,
$\psi_1(I)=\lambda I$, and $\psi_1=\lambda\psi_0$. Now if we apply
(13) we will obtain $I=\psi_0(I\cdot
I)=\psi_1(I)\psi_1(I)=\lambda^2 I$. In this case $\lambda=\pm 1$.
This argument allows us to conclude that for each automorphism
$\psi$ there is an automorphism $\psi_0$ of $B$ such that either
$\psi(X+tY)=\psi_0(X)+t\psi_0(Y)$ or
$\psi(X+tY)=\psi_0(X)-t\psi_0(Y)$. The proof is complete.
\end{proof}

\begin{lemma}\label{ltm22}
Let $A=B+tB$ where $B\cong M_n(F)$ be an associative superalgebra
of type $Q(n)$, and $\psi$ be an antiautomorphism of $A$. Then there
exists an antiautomorphism $\psi_0$ of $B$ such that for any $X+tY\in
A$, either $\psi(X+tY)=\psi_0(X)+it\psi_0(Y)$ or
$\psi(X+tY)=\psi_0(X)-it\psi_0(Y)$ where $i^2=-1$.
\end{lemma}

\begin{proof}
The proof of this lemma is similar to the previous one except that in the case where $\psi$
 is a (super!)antiautomorphism the equations (10--13) are replaced by

 $$\psi_0(X_1X_2)=\psi_0(X_2)\psi_0(X_1)\eqno(14)$$
 $$\psi_1(XY)=\psi_1(Y)\psi_0(X)\eqno(15)$$
 $$\psi_1(YX)=\psi_0(X)\psi_1(Y)\eqno(16)$$
 $$\psi_0(Y_1Y_2)=-\psi_1(Y_2)\psi_1(Y_1).\eqno(17)$$

where all $X_1,X_2,X, Y_1,Y_2,Y\in B$. Now (14) implies $\psi_0$
being an antiautomorphism. Also (15) and (16) imply
$\psi_1(I)=\lambda I$ and $\psi_1=\lambda\psi_0$. Using (17), we
now derive that $\lambda=\pm i$.
\end{proof}

Now we are ready to prove the second main result of this paper.

\begin{theorem}\label{tm2} Let $G$ be a finite abelian group and $R=A\oplus A^{sop}$ an involution simple superalgebra with $A$ of type $Q(n)$,
 as in item(3) of Proposition \ref{p1}. Suppose the base field $F$ is  algebraically closed of characteristic 0 or coprime to the order of $G$.Then
 any $G$-grading of   $A=B+tB\cong Q(n)$ with the
ordinary exchange involution $*$ compatible with $G$-grading has one
of the following forms:

\emph{Type I:} $R_g=A_g\oplus A_g^{sop}$, for  a grading of
$A=\oplus_{g\in G} A_g,$

\emph{Type II:} $R_g=\{(x,x^{\dagger})|\,x\in B_g\}\oplus
\{(tx,-tx^{\dagger})|\,x\in B_{gh}\}\oplus \{(x,-x^{\dagger})|\,x\in B_{gh^2}\}
\oplus \{(tx,tx^{\dagger})|\,x\in B_{gh^3}\}$, where $h$ is an element of order 4 in $G$, $\dagger$ is an involution on $B\cong M_n(F)$, $B=\bigoplus_{g\in G}B_g$ is an involution
grading on $B$ with respect to involution $\dagger$.
\end{theorem}

\begin{proof}

If $\wg$ acts on $R$ by the automorphisms of Type 1 only, we
arrive at Type I gradings described just before the statement of
Theorem \ref{tm1}. Otherwise, let $\alpha(\wg)$ contain all
possible automorphisms, where, as before,
$\alpha:\wg\rightarrow\aut{R}$ is the homomorphism corresponding
to our grading. Let $\vp\in\alpha(\wg)$ be such that
$\vp((x,y))=(\vp_0(y),\vp_0(x))$ where $\vp_0$ is an
antiautomorphism of $A$, $x,y\in A$. Then, according to Lemma
\ref{ltm22}, $\vp_0$ has one of two forms $\vp_0(u+tv)=\vp_1(u)\pm
it\vp_1(v)$, where $\vp_1$ is an antiatomorphism of $B=M_n(F)$,
and $x=u+tv$ with $u,v\in B$. If we compute the powers of $\vp$ on
$(x,y)=(u+tv,p+tq)$ where also $p,q\in B$ then we obtain the
following:

$$\vp^2((u+tv,p+tq))=(\vp_1^2(u)-t\vp^2_1(v),\vp_1^2(p)-t\vp^2_1(q))\eqno(18)$$
$$\vp^3((u+tv,p+tq))=(\vp_1^3(p)\mp t\vp^3_1(q),\vp_1^3(u)\mp
t\vp^3_1(v))\eqno(19)$$
$$\vp^4((u+tv,p+tq))=(\vp_1^4(u)+t\vp^4_1(v),\vp_1^4(p)+t\vp^4_1(q))\eqno(20)$$

Clearly, if we replace $\vp$ by $\vp^3$ we may assume from the
very beginning that $\vp_0(u+tv)=\vp_1(u)+ it\vp_1(v)$, for an
antiautomorphism $\vp_1$ of $B$. Let $\zeta\in\wg$ be such that
$\alpha(\zeta)=\vp$ and let
$$
\Lambda=\left\{\chi\in\wg\,|\,\alpha(\chi)((u+tv,p+tq))=(\pi_1(u)+t\pi_1(v),\pi_1(p)+t\pi_1(q))\right\}
$$
for any $u,v,p,q\in B$, $\pi_1\in\aut{B}$. Then
$\wg=\Lambda\cup\Lambda\zeta\cup\Lambda\zeta^2\cup\Lambda\zeta^3$.
Indeed, choose any $\chi\in\wg$ and consider $\pi=\alpha(\chi)$.
Then $\pi((x,y))$ is described by Lemmas \ref{l5} and then
\ref{ltm21} or \ref{ltm22}. Direct calculations using equations
(18,19,20) show that if $\pi$ is one of the cases of Lemma
\ref{ltm21} then either $\chi\in\Lambda$ or
$\vp^{2}\chi\in\Lambda$. If $\pi$ is one of the cases of Lemma
\ref{ltm22} then either $\vp\chi\in\Lambda$ or
$\vp^3\chi\in\Lambda$.

Let us define a mapping $\alpha_1:\wg\rightarrow\baut{B}$ by
associating with each $\chi\in\wg$ the mapping $\pi_1$ as in the
previous paragraph. Obviously, this is a homomorphism of groups
and the image $\vp_1$ of $\zeta$ is an antiautomorphism. In this
case Theorem \ref{tautant} applies and there exists an
automorphism $\psi_1$ of $B$ such that $\psi_1^2=\vp_1^2$ and
$\psi_1$ commutes with every $\pi_1\in\alpha_1(\wg)$. Let use our
previous notation to define an automorphism $\psi$ of $R$ by
setting $\psi((x,y))=(\psi_0(x),\psi_0(y))$ where
$\psi_0(u+tv)=\psi_1(u)+t\psi_1(v)$. Immediate calculations using
different cases of Lemmas \ref{ltm21} or \ref{ltm22} show that
$\psi$ commutes with any element of $\alpha(\wg)$. For example, if
$\pi\in\alpha(\wg)$ has the form
$\pi((u+tv,p+tq))=(\pi_1(p)-it\pi_1(q),\pi_1(u)-it\pi_1(v))$ then
using that $\psi_1\pi_1=\pi_1\psi_1$, we easily find both
$\psi\pi$ and $\pi\psi$ acting on $(u+tv,p+tq)$ produce the same
$(\psi_1\pi_1(p)- it\psi_1\pi_1(q),\psi_1\pi_1(u)-
it\psi_1\pi_1(v))$.

In order to apply Exchange Theorem, we define another mapping $\beta:\wg\rightarrow\aut{R}$ by
 setting $\beta(\zeta^k\lambda)=\psi^k\alpha(\lambda)$ for $k=0,1,2,3$. By Equation (19), $\vp^4=\psi^4$ and so this mapping is well defined and is a homomorphism of groups coinciding with $\alpha$ on $\Lambda$. Let also $R=\bigoplus_{g\in G}\overline{R}_g$ be a $G$-grading defined by $\beta$. This allows to apply Exchange Theorem, in which $\Lambda^{\perp}=H=\{e,h,h^2,h^3\}$. The homomorphism $\gamma:\wg\rightarrow\aut{R}$ defined by $\gamma(\chi)=\alpha^{-1}(\chi)\beta(\chi)$ defines a grading $R=R^{(e)}\oplus R^{(h)}\oplus R^{(h^2)}\oplus R^{(h^3)}$ and
$$
R_g=\overline{R}_g\cap R^{(e)}\oplus \overline{R}_{gh}\cap
R^{(h)}\oplus \overline{R}_{gh^2}\cap R^{(h^2)}\oplus
\overline{R}_{gh^3}\cap R^{(h^3)}.\eqno(21)
$$

Let us consider $\theta=\gamma(\zeta)$. Then the respective $\theta_1\in\baut{B}$ is an involution, which we denote by $\dagger$.
The grading of $R$ defined by $\beta$ induces a grading $B=\bigoplus_{g\in G}B_g$ on $B$, which permutes with $\dagger$,
 hence is an involution grading on the matrix algebra $B=M_n(F)$. We have $\overline{R}_g=\{(u+tv,p+tq)\,|\,u,v,p,q\in B_g\}$.
 To finally compute the homogeneous components of our original grading by Equation (21), we need to compute the components of the
 $H$-grading $R^{(t)}$, $t\in H$. We have $(x,y)\in R^{(e)}$ if $\theta((x,y))=(x,y)$. Using the same notation for $x,y\in A$, as before, we get
$$
\theta((u+tv,p+tq))=(\theta_1(p)-it\theta_1(q),\theta_1(u)-it\theta_1(v))=(p^{\dagger}-itq^{\dagger},u^{\dagger}-itv^{\dagger}).
$$
If $(x,y)\in R^{(e)}$ we must have $p^{\dagger}-itq^{\dagger}=u+tv$, $u^{\dagger}-itv^{\dagger}=p+tq$ and so $p^{\dagger}=u$, $-iq^{\dagger}=v$, $u^{\dagger}=p$, and $-iv^{\dagger}=q$. It follows that $p=u^{\dagger}$, $v=q=0$. Finally, $R^{(e)}=\{(u,u^{\dagger})\,|\, u\in B\}$. Now since $\overline{R}_g=\{(u+tv,p+tq)\,|\,u,v,p,q\in B_g\}$, we finally obtain $\overline{R}_g\cap R^{(e)}=\{(u,u^{\dagger})\,|\, u\in B\}$.

A similar computation gives us also $\overline{R}_{gh}\cap
R^{(h)}=\{(tv,-tv^{\dagger})\,|\, v\in B\}$,
$\overline{R}_{gh^2}\cap R^{(h^2)}=\{(u,-u^{\dagger})\,|\, u\in
B\}$, and $\overline{R}_{gh^3}\cap
R^{(h^3)}=\{(tv,tv^{\dagger})\,|\, v\in B\}$.

Now the proof of our theorem is complete.
\end{proof}

\end{document}